\documentclass[11pt]{amsart}

\newtheorem{lemma}{Lemma}
\newtheorem{theorem}{Theorem}
\newtheorem{proposition}{Proposition}

\theoremstyle{definition}

\makeindex

\begin{document}

\title[The sectional curvature of tangent bundles ...]
{THE SECTIONAL CURVATURE OF THE TANGENT BUNDLES WITH GENERAL
NATURAL LIFTED METRICS}
\author{S.~L.~Dru\c t\u a}

\begin{minipage}{2.5in}
\begin{flushleft}
Ninth International Conference on Geometry, Integrability and
Quantization, June 8–-13, 2007, Varna, Bulgaria \\Ivaïlo M.
Mladenov, Editor \\
SOFTEX, Sofia 2008, pp 198--209
\end{flushleft}
\end{minipage}

\begin{abstract}
We study some properties of the tangent bundles with metrics of
general natural lifted type. We consider a Riemannian manifold
$(M,g)$ and we find the conditions under which the Riemannian
manifold $(TM,G)$, where $TM$ is the tangent bundle of $M$ and $G$
is the general natural lifted metric of $g$, has constant
sectional curvature.

%\vskip 1cm
\rightskip=1.2cm \leftskip=1.2cm {\it Key words and phrases}:
tangent bundle, Riemannian metric, general natural lift, sectional
curvature.
\end{abstract}

\maketitle

\section{Introduction}\label{sec:1}

In the geometry of the tangent bundle $TM$ of a smooth
$n$-dimensional Riemannian manifold $(M,g)$ one uses several
Riemannian and pseudo-Rieman-\ nian metrics, induced by the
Riemannian metric $g$ on $M$. Among them, we may quote the Sasaki
metric, the Cheeger-Gromoll metric and the complete lift of the
metric $g$. The possibility to consider vertical, complete and
horizontal lifts on the tangent bundle $TM$ (see \cite{YanoIsh})
leads to some interesting geometric structures, studied in the
last years (see \cite{Abbassi1}, \cite{Abbassi2}, \cite{BejanOpr},
\cite{Munteanu1}, \cite{Munteanu2}, \cite{Tahara}), and to
interesting relations with some problems in Lagrangian and
Hamiltonian mechanics. On the other hand, the na-\ tural lifts of
$g$ to $TM$ (introduced in \cite{KowalskiSek} and \cite{Krupka})
induce some new Riemannian and pseudo-Riemannian geometric
structures with many nice geometric properties (see \cite{Kolar},
\cite{KowalskiSek}).

Professor Oproiu has studied some properties of a natural lift
$G$, of diagonal type, of the Riemannian metric $g$ and a natural
almost complex structure $J$ of diagonal type  on $TM$ (see
\cite{Oproiu2}, \cite{Oproiu3}, \cite{Oproiu4}, and see also
\cite{OprPap1}, \cite{OprPap2}). In the paper \cite{Oproiu1}, the
same author has presented a general expression of the natural
almost complex structures on $TM$. In the definition of the
natural almost complex structure $J$ of general type there are
involved eight parameters (smooth functions of the density energy
on $TM$). However, from the condition for $J$ to define an almost
complex structure, four of the above parameters can be expressed
as (rational) functions of the other four parameters. A Riemannian
metric $G$ which is a natural lift of general type of the metric
$g$ depends on other six parameters. In \cite{OprDruta} we have
found the conditions under which the K\"ahlerian manifold
$(TM,G,J)$ has constant holomorphic sectional curvature.

In the present paper we study the sectional curvature of the
tangent bundle of a Riemannian manifold $(M,g)$. Namely, we are
interested in finding the conditions under which the Riemannian
manifold $(TM,G)$, where $G$ is the general natural lifted metric
of $g$, has constant sectional curvature. We obtain that the
sectional curvature of $(TM,G)$ is zero and the base manifold must
be flat.
\section{Preliminary results}\label{sec:2}

Consider a smooth $n$-dimensional Riemannian manifold $(M,g)$ and
denote its tangent bundle by $\tau :TM\longrightarrow M$. Recall
that $TM$ has a structure of a $2n$-dimensional smooth manifold,
induced from the smooth manifold structure of $M$. This structure
is obtained by using local charts on $TM$ induced  from usual
local charts on $M$. If $(U,\varphi )= (U,x^1,\dots ,x^n)$ is a
local chart on $M$, then the corresponding induced local chart on
$TM$ is $(\tau ^{-1}(U),\Phi )=(\tau ^{-1}(U),x^1,\dots , x^n,$
$y^1,\dots ,y^n)$, where the local coordinates $x^i,y^j,\
i,j=1,\dots ,n$, are defined as follows. The first $n$ local
coordinates of a tangent vector $y\in \tau ^{-1}(U)$ are the local
coordinates in the local chart $(U,\varphi)$ of its base point,
i.e. $x^i=x^i\circ \tau$, by an abuse of notation. The last $n$
local coordinates $y^j,\ j=1,\dots ,n$, of $y\in \tau ^{-1}(U)$
are the vector space coordinates of $y$ with respect to the
natural basis in $T_{\tau(y)}M$ defined by the local chart
$(U,\varphi )$. Due to this special structure of differentiable
manifold for $TM$, it is possible to introduce the concept of
$M$-tensor field on it. The $M$-tensor fields are defined by their
components with respect to the induced local charts on $TM$ (hence
they are defined locally), but they can be interpreted as some
(partial) usual tensor fields on $TM$. However, the essential
quality of an $M$-tensor field on $TM$ is that the local
coordinate change rule of its components with respect to the
change of induced local charts is the same as the local coordinate
change rule of the components of an usual tensor field on $M$ with
respect to the change of local charts on $M$. More precisely, an
$M$-tensor field of type $(p,q)$ on $TM$ is defined by sets of
$n^{p+q}$ components (functions depending on $x^i$ and $y^i$),
with $p$ upper indices and $q$ lower indices, assigned to induced
local charts $(\tau ^{-1}(U),\Phi )$ on $TM$, such that the local
coordinate change rule of these components (with respect to
induced local charts on $TM$) is that of the local coordinate
components  of a tensor field of type $(p,q)$ on the base manifold
$M$ (with respect to usual local charts on $M$), when a change of
local charts on $M$ (and hence on $TM$) is performed (see
\cite{Mok} for further details); e.g., the components $y^i,\
i=1,\dots ,n$, corresponding to the last $n$ local coordinates of
a tangent vector $y$, assigned to the induced local chart $(\tau
^{-1}(U), \Phi )$ define an $M$-tensor field of type $(1,0)$ on
$TM$. An usual tensor field of type $(p,q)$ on $M$ may be thought
of as an $M$-tensor field of type $(p,q)$ on $TM$. If the
considered tensor field on $M$ is covariant only, the
corresponding $M$-tensor field on $TM$ may be identified with the
induced (pullback by $\tau $) tensor field on $TM$.  Some useful
$M$-tensor fields on $TM$ may be obtained as follows. Let
$u:[0,\infty )\longrightarrow {\mathbb R}$ be a smooth function
and let $\|y\|^2=g_{\tau (y)}(y,y)$ be the square of the norm of
the tangent vector $y\in \tau ^{-1}(U)$. If $\delta ^i_j$ are the
Kronecker symbols (in fact, they are the local coordinate
components of the identity tensor field $I$  on $M$), then the
components $u(\|y\|^2)\delta ^i_j$ define an $M$-tensor field of
type $(1,1)$ on $TM$. Similarly, if $g_{ij}(x)$ are the local
coordinate components of the metric tensor field $g$ on $M$ in the
local chart $(U,\varphi )$, then the components $u(\|y\|^2)
g_{ij}$ define a symmetric $M$-tensor field of type $(0,2)$ on
$TM$. The components $g_{0i}=y^kg_{ki}$ define an $M$-tensor field
of type $(0,1)$ on $TM$.

Denote by $\dot \nabla$ the Levi Civita connection of the
Riemannian metric $g$ on $M$. Then we have the direct sum
decomposition
\begin{equation}
TTM=VTM\oplus HTM
\end{equation}
of the tangent bundle to $TM$ into the vertical distribution
$VTM={\rm Ker}\ \tau_*$ and the horizontal distribution $HTM$
defined by $\dot \nabla $. The set of vector fields
$(\frac{\partial}{\partial y^1}, \dots , \frac{\partial}{\partial
y^n})$ on $\tau ^{-1}(U)$ defines a local frame field for $VTM$
and for $HTM$ we have the local frame field $(\frac{\delta}{\delta
x^1},\dots ,\frac{\delta}{\delta x^n})$, where
$$
\frac{\delta}{\delta x^i}=\frac{\partial}{\partial
x^i}-\Gamma^h_{0i} \frac{\partial}{\partial y^h},\ \ \ \Gamma
^h_{0i}=y^k\Gamma ^h_{ki}
 $$
and $\Gamma ^h_{ki}(x)$ are the Christoffel symbols of $g$.

The set $(\frac{\partial}{\partial y^1},\dots
,\frac{\partial}{\partial y^n}, \frac{\delta}{\delta x^1},\dots
,\frac{\delta}{\delta x^n})$ defines a local frame on $TM$,
adapted to the direct sum decomposition (1). Remark that
$$
\frac{\partial}{\partial y^i}=(\frac{\partial}{\partial x^i})^V,\
\ \frac{\delta}{\delta x^i}=(\frac{\partial}{\partial x^i})^H
$$
where $X^V$ and $X^H$ denote the vertical and horizontal lift of
the vector field $X$ on $M$ respectively. We can use the vertical
and horizontal lifts  in order to obtain invariant expressions for
some results in this paper. However, we should prefer to work in
local coordinates  since the formulas are obtained easier and, in
a certain sense, they are more natural.

We can easily obtain the following

\begin{lemma}\label{lema1}
If $n>1$ and $u,v$ are smooth functions on $TM$ such that
$$
u g_{ij}+v g_{0i}g_{0j}=0
$$
on the domain of any induced local chart on $TM$, then $u=0,\
v=0$.
\end{lemma}

\bf Remark. \rm In a similar way  we obtain from the condition
$$
u\delta ^i_j+vg_{0j}y^i=0
$$
the relation $u=v=0$.

Consider the energy density of the tangent vector $y$ with respect
to the Riemannian metric $g$
$$
t=\frac{1}{2}\|y\|^2=\frac{1}{2}g_{\tau(y)}(y,y)=\frac{1}{2}g_{ik}(x)y^iy^k,
\ \ \ y\in \tau^{-1}(U).
$$
Obviously, we have $t\in [0,\infty)$ for all $y\in TM$.

\section{The sectional curvature of the tangent bundle with general natural lifted metric}

Let $G$ be the general natural lifted metric on $TM$, defined by
\begin{equation}\label{defG}
\begin{array}{l}
G(\frac{\delta}{\delta x^i}, \frac{\delta}{\delta x^j})=c_1g_{ij}+
d_1g_{0i}g_{0j}=G^{(1)}_{ij}
\\   \mbox{ } \\
G(\frac{\partial}{\partial y^i}, \frac{\partial}{\partial y^j})=
c_2g_{ij}+d_2g_{0i}g_{0j}=G^{(2)}_{ij}
\\   \mbox{ } \\
G(\frac{\partial}{\partial y^i},\frac{\delta}{\delta x^j})=
G(\frac{\delta}{\delta x^i},\frac{\partial}{\partial y^j})=
c_3g_{ij}+d_3g_{0i}g_{0j}=G^{(3)}_{ij}
\end{array}
\end{equation}
where $c_1,c_2,c_3,d_1,d_2,d_3$ are six smooth functions of the
density energy on $TM$.

The Levi-Civita connection $\nabla$ of the Riemannian manifold
$(TM,G)$ is obtained from the formula
$$
2G(\nabla_XY,Z)=X(G(X,Z))+Y(G(X,Z))-Z(G(X,Y))
$$
$$
+G([X,Y],Z)-G([X,Z],Y)-G([Y,Z],X); ~~\forall X,Y,Z\in \chi(M)
$$
and is characterized by the conditions
$$
\nabla G=0,\ T=0
$$
where $T$ is the torsion tensor of $\nabla.$

In the case of the tangent bundle $TM$ we can obtain the explicit
expression of $\nabla$. The symmetric $2n\times 2n$ matrix
$$
\left(
\begin{array}{ll}
G^{(1)}_{ij} &  G^{(3)}_{ij} \\
G^{(3)}_{ij} &  G^{(2)}_{ij}
\end{array}
\right)
$$
associated to the metric $G$ in the base $(\frac{\delta}{\delta
x^1},\dots ,\frac{\delta}{\delta x^n},\frac{\partial}{\partial
y^1},\dots ,\frac{\partial}{\partial y^n})$ has the inverse
$$
\left(
\begin{array}{ll}
H_{(1)}^{ij} &  H_{(3)}^{ij} \\
H_{(3)}^{ij} &  H_{(2)}^{ij}
\end{array}
\right)
$$
where the entries are the blocks
 \begin{eqnarray*}
H_{(1)}^{kl}=p_1g^{kl}+q_1y^ky^l
\end{eqnarray*}
\begin{eqnarray}\label{matrinv}
H_{(2)}^{kl}=p_2g^{kl}+q_2y^ky^l
\end{eqnarray}
\begin{eqnarray*}
H_{(3)}^{kl}=p_3g^{kl}+q_3y^ky^l.
\end{eqnarray*}
Here $g^{kl}$ are the components of the inverse of the matrix
$(g_{ij})$ and $p_1,q_1,p_2,q_2,p_3$, $q_3:[0,\infty)\rightarrow
\mathbb R,$ some real smooth functions. Their expressions are
obtained  by solving the system:
$$
\begin{cases}
G^{(1)}_{ih}H_{(1)}^{hk}+G^{(3)}_{ih}H_{(3)}^{hk}=\delta_i^k\\
G^{(1)}_{ih}H_{(3)}^{hk}+G^{(3)}_{ih}H_{(2)}^{hk}=0\\
G^{(3)}_{ih}H_{(1)}^{hk}+G^{(2)}_{ih}H_{(3)}^{hk}=0\\
G^{(3)}_{ih}H_{(3)}^{hk}+G^{(2)}_{ih}H_{(2)}^{hk}=\delta_i^k
\end{cases}
$$
in which we substitute the relations (\ref{defG}) and
(\ref{matrinv}). By using Lemma \ref{lema1},  we get $p_1,p_2,p_3$
as functions of $c_1,c_2,c_3$

\begin{eqnarray}\label{inversa1}
p_1=\frac{c_2}{c_1c_2-c_3^2},\ \ p_2=\frac{c_1}{c_1c_2-c_3^2},\ \
p_3=-\frac{c_3}{c_1c_2-c_3^2}
\end{eqnarray}
and $q_1,q_2,q_3$ as functions of $c_1,c_2,c_3,$ $d_1,d_2,d_3,$
$p_1,p_2,p_3$
$$
q_1=-\frac{c_2d_1p_1 - c_3d_3p_1 - c_3d_2p_3 + c_2d_3p_3 +
2d_1d_2p_1t - 2d_3^2p_1t}{c_1c_2 - c_3^2 + 2c_2d_1t + 2c_1d_2t -
4c_3d_3t + 4d_1d_2t^2 - 4d_3^2t^2}
$$
\begin{eqnarray}\label{inversa2}
q_2=-\frac{d_2p_2 + d_3p_3}{c_2 + 2d_2t}
\end{eqnarray}
$$
+\frac{(c_3 + 2d_3t)[(d_3p_1 + d_2p_3)(c_1 + 2d_1t) - (d_1p_1 +
d_3p_3)(c_3 + 2d_3t)]}{(c_2 + 2d_2t)[(c_1 + 2d_1t)(c_2 + 2d_2t) -
(c_3 + 2d_3t)^2]}
$$
$$
q_3=-\frac{(d_3p_1 + d_2p_3)(c_1 + 2d_1t) - (d_1p1 + d_3p_3)(c_3 +
2d_3t)}{(c_1 + 2d_1t)(c_2 + 2d_2t) - (c_3 + 2d_3t)^2}.
$$
In the paper \cite{OprDruta} we obtained the expression of the
Levi Civita connection of the Riemannian metric $G$ on $TM$.
\begin{theorem}
The Levi-Civita connection $\nabla$ of\ $G$ has the following
expression in the local adapted frame $(\frac{\partial}{\partial
y^1}, \dots , \frac{\partial}{\partial y^n},\frac{\delta}{\delta
x^1},\dots ,\frac{\delta}{\delta x^n})$
$$
\begin{cases}
$$
\nabla_{\frac{\partial}{\partial y^i}} \frac{\partial}{\partial
y^j}=Q^{h}_{ij}\frac{\partial}{\partial
y^h}+\widetilde{Q}^{h}_{ij}\frac{\delta}{\delta x^h},~
\nabla_{\frac{\delta}{\delta x^i}} \frac{\partial}{\partial
y^j}=(\Gamma^h_{ij}+\widetilde{P}^{h}_{ji})\frac{\partial}{\partial
y^h}+P^{h}_{ji}\frac{\delta}{\delta x^h}\\

\nabla_{\frac{\partial}{\partial y^i}} \frac{\delta}{\delta
x^j}=P^{h}_{ij}\frac{\delta}{\delta
x^h}+\widetilde{P}^{h}_{ij}\frac{\partial}{\partial y^h},~
\nabla_{\frac{\delta}{\delta x^i}} \frac{\delta}{\delta
x^j}=(\Gamma^h_{ij}+\widetilde{S}^{h}_{ij})\frac{\delta}{\delta
x^h}+S^{h}_{ij}\frac{\partial}{\partial y^h}
$$
\end{cases}
$$
where $\Gamma^h_{ij}$ are the Christoffel symbols of the
connection $\dot\nabla$ and the $M$-tensor fields appearing as
coefficients in the above expressions are given as
$$
\left\{
\begin{array}{l}
Q_{ij}^h=\frac{1}{2}(\partial_iG_{jk}^{(2)}+\partial_jG_{ik}^{(2)}-
\partial_kG_{ij}^{(2)})H_{(2)}^{kh}+\frac{1}{2}(\partial_iG_{jk}^{(3)}+
\partial_jG_{ik}^{(3)})H_{(3)}^{kh}\\
\widetilde{Q}_{ij}^h=\frac{1}{2}(\partial_iG_{jk}^{(2)}+\partial_jG_{ik}^{(2)}-
\partial_kG_{ij}^{(2)})H_{(3)}^{kh}+\frac{1}{2}(\partial_iG_{jk}^{(3)}+
\partial_jG_{ik}^{(3)})H_{(1)}^{kh}\\\\
P^h_{ij}=\frac{1}{2}(\partial_iG_{jk}^{(3)}-
\partial_kG_{ij}^{(3)})H_{(3)}^{kh}+\frac{1}{2}(\partial_iG_{jk}^{(1)}+
R^l_{0jk}G^{(2)}_{li})H_{(1)}^{kh}\\
\widetilde{P}^h_{ij}=\frac{1}{2}(\partial_iG_{jk}^{(3)}-
\partial_kG_{ij}^{(3)})H_{(2)}^{kh}+\frac{1}{2}(\partial_iG_{jk}^{(1)}+
R^l_{0jk}G^{(2)}_{li})H_{(3)}^{kh}\\\\
S^h_{ij}=-\frac{1}{2}(\partial_kG_{ij}^{(2)}+R^l_{0ij}G^{(2)}_{lk})H_{(2)}^{kh}+c_3R_{i0jk}H_{(3)}^{kh}\\
\widetilde{S}^h_{ij}=-\frac{1}{2}(\partial_kG_{ij}^{(1)}+R^l_{0ij}G^{(2)}_{lk})H_{(3)}^{kh}+c_3R_{i0jk}H_{(1)}^{kh}
\end{array}
\right.
$$
where $R^h_{kij}$ are the components of the curvature tensor field
of the Levi Civita connection $\dot \nabla$ of the base manifold
$(M,g)$.
\end{theorem}

Taking into account the expressions (\ref{defG}), (\ref{matrinv})
 and by using the formulas (\ref{inversa1}), (\ref{inversa2}) we
can obtain the detailed expressions of
$P^h_{ij},Q^h_{ij},S^h_{ij}, \widetilde P^h_{ij}, \widetilde
Q^h_{ij}, \widetilde S^h_{ij}.$

The curvature tensor field $K$ of the connection $\nabla$ is
defined by the well known formula
$$
K(X,Y)Z=\nabla_X\nabla_YZ-\nabla_Y\nabla_XZ-\nabla_{[X,Y]}Z,\ \ \
X,Y,Z\in \Gamma (TM).
$$

By using the local adapted frame $(\frac{\delta}{\delta
x^i},\frac{\partial}{\partial y^j}),\ i,j=1,\dots ,n$ we obtained
in \cite{OprDruta}, after a standard straightforward computation
$$
K\big(\frac{\delta}{\delta x^i},\frac{\delta}{\delta x^j}
\big)\frac{\delta}{\delta x^k}=XXXX^h_{kij}\frac{\delta}{\delta
x^h}+XXXY^h_{kij}\frac{\partial}{\partial y^h}
$$
$$
K\big(\frac{\delta}{\delta x^i},\frac{\delta}{\delta x^j}
\big)\frac{\partial}{\partial
y^k}=XXYX^h_{kij}\frac{\delta}{\delta
x^h}+XXYY^h_{kij}\frac{\partial}{\partial y^h}
$$
$$
K\big(\frac{\partial}{\partial y^i},\frac{\partial}{\partial y^j}
\big)\frac{\delta}{\delta x^k}=YYXX^h_{kij}\frac{\delta}{\delta
x^h}+YYXY^h_{kij}\frac{\partial}{\partial y^h}
$$
$$
K\big(\frac{\partial}{\partial y^i},\frac{\partial}{\partial y^j}
\big)\frac{\partial}{\partial
y^k}=YYYX^h_{kij}\frac{\delta}{\delta
x^h}+YYYY^h_{kij}\frac{\partial}{\partial y^h}
$$
$$
K\big(\frac{\partial}{\partial y^i},\frac{\delta}{\delta x^j}
\big)\frac{\delta}{\delta x^k}=YXXX^h_{kij}\frac{\delta}{\delta
x^h}+YXXY^h_{kij}\frac{\partial}{\partial y^h}
$$
$$
K\big(\frac{\partial}{\partial y^i},\frac{\delta}{\delta x^j}
\big)\frac{\partial}{\partial
y^k}=YXYX^h_{kij}\frac{\delta}{\delta
x^h}+YXYY^h_{kij}\frac{\partial}{\partial y^h}
$$\vskip1mm
where the $M$-tensor fields appearing as coefficients denote the
horizontal and vertical components of the curvature tensor of the
tangent bundle, and they are given by
$$
XXXX^h_{kij}=\widetilde{S}^h_{il}\widetilde{S}^l_{jk}+P^h_{li}S^l_{jk}
-\widetilde{S}^h_{jl}\widetilde{S}^l_{ik}-P^h_{lj}S^l_{ik}+R^h_{kij}+R_{0ij}^lP^h_{lk}
$$
$$
XXXY^h_{kij}=\widetilde{S}^l_{jk}S^h_{il}+\widetilde{P}^h_{li}S^l_{jk}-
\widetilde{S}^l_{ik}S^h_{jl}-\widetilde{P}^h_{lj}S^l_{ik}+\widetilde{P}^h_{lk}R^l_{0ij}
$$
$$
-\frac{1}{2}\dot{\nabla}_iR_{0jk}^rG^{(2)}_{rl}H^{(3)}_{hl}+c_3\dot{\nabla}_iR_{j0kh}
$$
$$
XXYX^h_{kij}=\widetilde{P}^l_{kj}P^h_{li}+P^l_{kj}\widetilde{S}^h_{il}-
\widetilde{P}^l_{ki}P^h_{lj}-P^l_{ki}\widetilde{S}^h_{jl}+R^l_{0ij}\widetilde{Q}^h_{lk}
$$
$$
XXYY^h_{kij}=\widetilde{P}^l_{kj}\widetilde{P}^h_{li}+P^l_{kj}S^h_{il}-
\widetilde{P}^l_{ki}\widetilde{P}^h_{lj}-P^l_{ki}S^h_{jl}+R^l_{0ij}Q^h_{lk}+R^h_{kij}
$$
$$
YYXX^h_{kij}=\partial_iP^h_{jk}-\partial_jP^h_{ik}+\widetilde{P}^l_{jk}\widetilde{Q}^h_{il}+
P^l_{jk}P^h_{il}-\widetilde{P}^l_{ik}\widetilde{Q}^h_{jl}-P^l_{ik}P^h_{jl}
$$
$$
YYXY^h_{kij}=\partial_i\widetilde{P}^h_{jk}-\partial_j\widetilde{P}^h_{ik}+\widetilde{P}^l_{jk}Q^h_{il}
+P^l_{jk}\widetilde{P}^h_{il}-\widetilde{P}^l_{ik}Q^h_{jl}-P^l_{ik}\widetilde{P}^h_{jl}
$$
$$
YYYX^h_{kij}=\partial_i\widetilde{Q}^h_{jk}-\partial_j\widetilde{Q}^h_{ik}+
Q^l_{jk}\widetilde{Q}^h_{il}+\widetilde{Q}^l_{jk}P^h_{il}-Q^l_{ik}\widetilde{Q}^h_{jl}
-\widetilde{Q}^l_{ik}P^h_{jl}
$$
$$
YYYY^h_{kij}=\partial_iQ^h_{jk}-\partial_jQ^h_{ik}+Q^l_{jk}Q^h_{il}+\widetilde{Q}^l_{jk}\widetilde{P}^h_{il}-
Q^l_{ik}Q^h_{jl}-\widetilde{Q}^l_{ik}\widetilde{P}^h_{jl}
$$
$$
YXXX^h_{kij}=\partial_i\widetilde{S}^h_{jk}+S^l_{jk}\widetilde{Q}^h_{il}+\widetilde{S}^l_{jk}P^h_{il}
-\widetilde{P}^l_{ik}P^h_{lj}-P^l_{ik}\widetilde{S}^h_{jl}-\dot{\nabla}_jR_{0ik}^rG^{(2)}_{rl}H^{(3)}_{hl}
$$
$$
YXXY^h_{kij}=\partial_iS^h_{jk}+S^l_{jk}Q^h_{il}+\widetilde{S}^l_{jk}\widetilde{P}^h_{il}
-\widetilde{P}^l_{ik}\widetilde{P}^h_{lj}-P^l_{ik}S^h_{jl}-\dot{\nabla}_jR_{0ik}^rG^{(2)}_{rl}H^{(1)}_{hl}
$$
$$
YXYX^h_{kij}=\partial_iP^h_{kj}+\widetilde{P}^l_{kj}\widetilde{Q}^h_{il}+P^l_{kj}P^h_{il}-
Q^l_{ik}P^h_{lj}-\widetilde{Q}^l_{ik}\widetilde{S}^h_{jl}
$$
$$
YXYY^h_{kij}=\partial_i\widetilde{P}^h_{kj}+\widetilde{P}^l_{kj}Q^h_{il}+P^l_{kj}\widetilde{P}^h_{il}-
Q^l_{ik}\widetilde{P}^h_{lj}-\widetilde{Q}^l_{ik}S^h_{jl}.
$$
We mention that we used the character $X$ on a certain position to
indicate that the argument on that position was a horizontal
vector field and, similarly, we used the character $Y$ for
vertical vector fields. \vskip1mm We compute the partial
derivatives with respect to the tangential coordinates $y^i$ of of
$G^{(\alpha)}_{jk}$ and $H_{(\alpha)}^{jk}$, for $\alpha =1,2,3$.
$$
\partial_iG^{(\alpha)}_{jk}=c_\alpha'g_{0i}g_{jk}+d_\alpha'g_{0i}g_{0j}g_{0k}+
d_\alpha g_{ij}g_{0k}+d_\alpha g_{0i}g_{jk}
$$

$$
\partial_iH_{(\alpha)}^{jk}=p_\alpha'g^{jk}g_{0i}+q_\alpha'g_{0i}y^jy^k+
q_\alpha\delta^j_iy^k+q_\alpha y^j\delta^k_i
$$

\vskip1mm
$$
\partial_i\partial_jG^{(\alpha)}_{kl}=c_\alpha''g_{0i}g_{0j}g_{kl}+c_\alpha'g_{ij}g_{kl}+
d_\alpha''g_{0j}g_{0k}g_{0l}+ d_\alpha'g_{ij}g_{0k}g_{0l}
$$

$$
+d_\alpha'g_{0j}g_{ik}g_{0l}+
d_\alpha'g_{0j}g_{0k}g_{il}+d_\alpha'g_{0i}g_{jk}g_{0l}+d_\alpha'g_{0i}g_{0k}g_{jl}+
d_\alpha g_{jk}g_{il}+d_\alpha g_{ik}g_{jl}.
$$
\vskip1mm Next we get the first order partial derivatives with
respect to the tangential coordinates $y^i$ of the $M$-tensor
fields $P^h_{ij},Q^h_{ij},S^h_{ij}, \widetilde P^h_{ij},
\widetilde Q^h_{ij}, \widetilde S^h_{ij}$
$$
\partial_iQ^h_{jk}=\frac{1}{2}\partial_iH_{(2)}^{hl}(\partial_jG^{(2)}_{kl}+\partial_kG^{(2)}_{jl}-
\partial_lG^{(2)}_{jk})+
\frac{1}{2}H_{(2)}^{hl}(\partial_i\partial_jG^{(2)}_{kl}+\partial_i\partial_kG^{(2)}_{jl}
$$
$$
-\partial_i\partial_lG^{(2)}_{jk})+\frac{1}{2}\partial_iH_{(3)}^{hl}(\partial_jG^{(3)}_{kl}+\partial_kG^{(3)}_{jl})+
\frac{1}{2}H_{(3)}^{hl}(\partial_i\partial_jG^{(3)}_{kl}+\partial_i\partial_kG^{(3)}_{jl})
$$
$$
\partial_i\widetilde{Q}^h_{jk}=\frac{1}{2}\partial_iH_{(3)}^{hl}(\partial_jG^{(2)}_{kl}+\partial_kG^{(2)}_{jl}-
\partial_lG^{(2)}_{jk})+
\frac{1}{2}H_{(3)}^{hl}(\partial_i\partial_jG^{(2)}_{kl}+\partial_i\partial_kG^{(2)}_{jl}
$$
$$
-\partial_i\partial_lG^{(2)}_{jk})+\frac{1}{2}\partial_iH_{(1)}^{hl}(\partial_jG^{(3)}_{kl}+\partial_kG^{(3)}_{jl})+
\frac{1}{2}H_{(1)}^{hl}(\partial_i\partial_jG^{(3)}_{kl}+\partial_i\partial_kG^{(3)}_{jl})
$$
$$
\partial_i\widetilde{P}^h_{jk}=\frac{1}{2}\partial_iH_{(2)}^{hl}(\partial_jG^{(3)}_{kl}-
\partial_lG^{(3)}_{jk})+
\frac{1}{2}H_{(2)}^{hl}(\partial_i\partial_jG^{(3)}_{kl}
-\partial_i\partial_lG^{(3)}_{jk})
$$
$$
+\frac{1}{2}\partial_iH_{(3)}^{hl}(\partial_jG^{(1)}_{kl}+R^r_{0kl}G^{(2)}_{rj})+
\frac{1}{2}H_{(3)}^{hl}(\partial_i\partial_jG^{(1)}_{kl}+R^r_{ikl}G^{(2)}_{rj}+R_{0kl}^r\partial_iG^{(2)}_{rj})
$$
$$
\partial_iP^h_{jk}=\frac{1}{2}\partial_iH_{(3)}^{hl}(\partial_jG^{(3)}_{kl}-
\partial_lG^{(3)}_{jk})+
\frac{1}{2}H_{(3)}^{hl}(\partial_i\partial_jG^{(3)}_{kl}
-\partial_i\partial_lG^{(3)}_{jk})
$$
$$
+\frac{1}{2}\partial_iH_{(1)}^{hl}(\partial_jG^{(1)}_{kl}+R^r_{0kl}G^{(2)}_{rj})+
\frac{1}{2}H_{(1)}^{hl}(\partial_i\partial_jG^{(1)}_{kl}+R^r_{ikl}G^{(2)}_{rj}+R_{0kl}^r\partial_iG^{(2)}_{rj})
$$
$$
\partial_iS^h_{jk}=-\frac{1}{2}\{(\partial_i\partial_rG^{(1)}_{jk}+R^l_{ijk}G^{(2)}_{lr}+
R_{0jk}^l\partial_iG^{(2)}_{lr})H_{(2)}^{rh}+
$$
$$
+(\partial_rG^{(1)}_{jk}+R_{0jk}^lG^{(2)}_{lr})\partial_iH_{(2)}^{rh}\}
+c_3'g_{0i}R_{j0kr}H_{(3)}^{rh}+c_3(R_{jikr}H_{(3)}^{rh}+R_{j0kr}\partial_iH_{(3)}^{rh})
$$
$$
\partial_i\widetilde{S}^h_{jk}=-\frac{1}{2}\{(\partial_i\partial_rG^{(1)}_{jk}+R^l_{ijk}G^{(2)}_{lr}
+R_{0jk}^l\partial_iG^{(2)}_{lr})H_{(3)}^{rh}+
$$
$$
+(\partial_rG^{(1)}_{jk}+R_{0jk}^lG^{(2)}_{lr})\partial_iH_{(3)}^{rh}\}
+c_3'g_{0i}R_{j0kr}H_{(1)}^{rh}+c_3(R_{jikr}H_{(1)}^{rh}+R_{j0kr}\partial_iH_{(1)}^{rh}).
$$

It was not convenient to think $c_1,c_2,c_3,d_1,d_2,d_3$ and
$p_1,p_2,p_3$, $q_1,q_2,q_3$ as functions of $t$ since RICCI did
not make some useful factorizations after the command
TensorSimplify. We decided to consider these functions as well as
their derivatives of first, second and  third order, as constants,
the tangent vector $y$ as a first order tensor, the components
$G^{(1)}_{ij}, G^{(2)}_{ij}, G^{(3)}_{ij},$ $H_{(1)}^{ij},
H_{(2)}^{ij}, H_{(3)}^{ij}$ as second order tensors and so on, on
the  Riemannian manifold $M$, the associated indices being
$h,i,j,k,l,r,s.$

\vskip5mm The  tensor field corresponding to the curvature tensor
field of a Riemannian manifold $(TM,G)$ having constant sectional
curvature $k$, is given by the formula:
$$
K_0(X,Y)Z=k[G(Y,Z)X-G(X,Z)Y].
$$

After a standard straightforward computation we obtain
$$
K_0\big(\frac{\delta}{\delta x^i},\frac{\delta}{\delta x^j}
\big)\frac{\delta}{\delta x^k}=XXXX0^h_{kij}\frac{\delta}{\delta
x^h}+XXXY0^h_{kij}\frac{\partial}{\partial y^h}
$$
$$
K_0\big(\frac{\delta}{\delta x^i},\frac{\delta}{\delta x^j}
\big)\frac{\partial}{\partial
y^k}=XXYX0^h_{kij}\frac{\delta}{\delta
x^h}+XXYY0^h_{kij}\frac{\partial}{\partial y^h}
$$
$$
K_0\big(\frac{\partial}{\partial y^i},\frac{\partial}{\partial
y^j} \big)\frac{\delta}{\delta
x^k}=YYXX0^h_{kij}\frac{\delta}{\delta
x^h}+YYXY0^h_{kij}\frac{\partial}{\partial y^h}
$$
$$
K_0\big(\frac{\partial}{\partial y^i},\frac{\partial}{\partial
y^j} \big)\frac{\partial}{\partial
y^k}=YYYX0^h_{kij}\frac{\delta}{\delta
x^h}+YYYY0^h_{kij}\frac{\partial}{\partial y^h}
$$
$$
K_0\big(\frac{\partial}{\partial y^i},\frac{\delta}{\delta x^j}
\big)\frac{\delta}{\delta x^k}=YXXX0^h_{kij}\frac{\delta}{\delta
x^h}+YXXY0^h_{kij}\frac{\partial}{\partial y^h}
$$
$$
K_0\big(\frac{\partial}{\partial y^i},\frac{\delta}{\delta x^j}
\big)\frac{\partial}{\partial
y^k}=YXYX0^h_{kij}\frac{\delta}{\delta
x^h}+YXYY0^h_{kij}\frac{\partial}{\partial y^h}
$$
where the $M$-tensor fields appearing as coefficients are the
horizontal and vertical components of the tensor $K_0$ and they
are given by
$$
XXXX0^h_{kij}=k[G^{(1)}_{jk}\delta_i^h-G^{(1)}_{ik}\delta^h_j],\quad
XXXY0^h_{kij}=0
$$
$$
XXYX0^h_{kij}=k[G^{(3)}_{jk}\delta_i^h-G^{(3)}_{ik}\delta^h_j],\quad
XXYY0^h_{kij}=0
$$
$$
YYXX0^h_{kij}=0,\quad
YYXY0^h_{kij}=k[G^{(3)}_{jk}\delta_i^h-G^{(3)}_{ik}\delta^h_j]
$$
$$
YYYX0^h_{kij}=0,\quad YYYY0^h_{kij}=k
[G^{(2)}_{jk}\delta_i^h-G^{(2)}_{ik}\delta^h_j]
$$
$$
YXXX0^h_{kij}=-kG^{(3)}_{ik}\delta_j^h,\quad
YXXY0^h_{kij}=kG^{(1)}_{jk}\delta_i^h
$$
$$
YXYX0^h_{kij}=-kG^{(2)}_{ik}\delta^h_j,\quad
YXYY0^h_{kij}=kG^{(3)}_{jk}\delta_i^h.
$$
\vskip2mm

In order to get the conditions under which $(TM,G)$ is a
Riemannian manifold of constant sectional curvature, we study the
vanishing of the components of the difference $K-K_0$. In this
study it is  useful the following generic result similar to the
lemma \ref{lema1}.

\begin{lemma}\label{lema2}
If $\alpha _1,\dots , \alpha_{10}$ are smooth functions on $TM$
such that
$$
\alpha_1 \delta^h_i g_{jk}+\alpha_2 \delta^h_j g_{ik}+ \alpha_3
\delta^h_kg_{ij}+\alpha_4 \delta^h_k g_{0i}g_{0j} +\alpha_5
\delta^h_j g_{0i}g_{0k}+\alpha_6 \delta^h_i g_{0j}g_{0k}
$$
$$
+\alpha_7g_{jk} g_{0i}y^h+ \alpha_8 g_{ik}g_{0j}y^h+\alpha_9
g_{ij}g_{0k}y^h+\alpha_{10}g_{0i}g_{0j}g_{0k}y^h=0
$$
then $\alpha_1=\dots =\alpha_{10}=0$.
\end{lemma}

After a detailed analysis of several terms in the vanishing
problem of the components of the above difference we can formulate
the next proposition.

\begin{proposition}\label{vanish of R}
Let $(M,g)$ be a Riemannian manifold. If the tangent bundle $TM$
with the general natural lifted metric $G$ has constant sectional
curvature, then the base manifold is flat.
\end{proposition}
\emph{Proof:} For $y=0$ the difference
$XXYY^{h}_{kij}-XXYY0^{h}_{kij}$ reduces to $R^h_{kij}.$

If the sectional curvature of the tangent bundle is constant, this
difference vanishes, so the curvature of the base manifold must
vanish too.

\section{Tangent bundles with constant sectional curvature}\label{sec:3}

\begin{theorem}\label{teorfinala}
Let $(M,g)$ be a Riemannian manifold. The tangent bundle $TM$ with
the natural lifted metric $G$ has constant sectional curvature if
and only if the base manifold is flat and the metric $G$ has the
associated matrix of the form:
$$
\begin{pmatrix}
c g_{ij} & \beta g_{ij}+
\beta' g_{0i} g_{0j} \\
\beta g_{ij}+ \beta' g_{0i} g_{0j} &\alpha
g_{ij}+\frac{\alpha'\beta^2+2\alpha'\beta
\beta't-2\alpha\beta'^2t}{\beta^2} g_{0i} g_{0j}
\end{pmatrix}
$$
where $\alpha,\beta$ are two real smooth function depending on the
energy density and $c$ is an arbitrary constant. Moreover, in this
case, $TM$ is flat, i.e. $k=0$.
\end{theorem}

\emph{Proof:} In the proposition \ref{vanish of R} we prooved that
the base manifold of the tangent bundle with constant sectional
curvature must be flat. By using the RICCI package of the program
Mathematica, we impose the vanishing condition for the curvature
tensor of the base maniflod in all the differences between the
components of the curvature tensors $K$ and $K_0$ of $TM$. After a
long computation we find some differences in which the third terms
are of one of the forms: $\frac{c_3d_1}{2
(c_3^2-c_1c_2))}g_{ij}\delta^h_k$ in the case of the differences
$YXXX^{h}_{kij}-YXXX0^{h}_{kij}$ and
$YXYY^{h}_{kij}-YXYY0^{h}_{kij}$, $\frac{c_1d_1}{2(c_1c_2 -
c_3^2)}g_{ij}\delta^h_k$ for the difference
$YXXY^{h}_{kij}-YXXY0^{h}_{kij}$ and $\frac{c_2d_1}{2(c_1c_2 -
c_3^2)}g_{ij}\delta^h_k$ for $YXYX^h_{kij}-YXYX^h_{kij}.$

As all the coeficients which appear in these differences must
vanish, we obtain $d_1=0$, because $c_1$ and $c_3$, or $c_2$ and
$c_3$ cannot vanish at the same time, the metric $g$ being
non-degenerated.

If we impose $d_1=0$ in $XXXY^{h}_{kij}-XXXY0^{h}_{kij}$ we obtain
that this difference contains the factors
$c_1c_1'(c_1'c_3-c_1c_3'+c_1d_3).$ Thus, for the annulation of
this difference, we have the cases
$c_1'=\frac{c_1c_3'-c_1d_3}{c_3}$ or $c_1=constant$ ($c_1=0$ being
a particular case).

The first case, $c_1'=\frac{c_1c_3'-c_1d_3}{c_3}$ is not a
favourable one, because the difference
$YYYY^{h}_{kij}-YYYY0^{h}_{kij}$ containes two summands which
cannot vanish:
$$
\frac{1}{2t}g_{jk}\delta^h_i - \frac{1}{2t}g_{ik}\delta^h_j.
$$
In the case $c_1=constant$ we obtain
$$
XXXX^{h}_{kij}-XXXX0^{h}_{kij}=-c_1k(g_{jk}\delta^h_i-g_{ik}\delta^h_j)
$$
from which $c_1=0$ or $k=0$. If $c_1=0$
$$
XXYX^{h}_{kij}-XXYX0^{h}_{kij}=-k(c_3g_{jk}\delta^h_i - c_3
g_{ik}\delta^h_j - d_3\delta^h_jg_{0i}g_{0k} +
d_3\delta^h_ig_{0j}g_{0k}).
$$
As we considered $c_1=0,$ we cannot have $c_3=0$ because the
metric $g$ must be non-degenerated, so the parenthesis cannot
vanish and it remaines $k=0.$ Now we can conclude that \emph{the
tangent bundle with general natural lifted metric cannot have
nonzero sectional curvature}.

We continue the study of the general case $c_1=constant$, since
the case $c_1=0$ is a particular case  only. Because the sectional
curvature of the tangent bundle, $k$, is null, we obtain that the
difference $XXYY^{h}_{kij}-XXYY0^{h}_{kij}$ vanishes if and only
if $d_3=c_3'.$ This condition makes vanish all the differences
that we study, except $YYYX^{h}_{kij}-YYYX0^{h}_{kij}.$ From the
annulation of this last difference, we obtain
$$
d_2=c_2'+2t\frac{c_2'c_3c_3'-c_2c_3'^2}{c_3^2}.
$$

If we denote $c_1$ by $c$, $c_2$ by $\alpha$ and $c_3$ by $\beta$,
we obtain that the matrix associated to the metric $G$ has the
form given in the theorem \ref{teorfinala}.

Therefore, the theorem \ref{teorfinala} gives the unique form of
the matrix associated to the metric $G$.

\section*{Acknowledgements}

Partially supported by the Grant AT No.191/2006, CNCSIS,
Ministerul Educa\c tiei \c si Cercet\u arii, Rom\^ania.

The author expresses her gratitude to the organizers, especially
to Professor Mladenov, for the Conference Grant.

S.L. Dru\c t\u a

Faculty of Mathematics

"Al.I. Cuza" University of Ia\c si

Bd. Carol I, no. 11

700506 Ia\c si

ROMANIA

\email{simonadruta@yahoo.com}

\end{document}